\documentclass[11pt]{article}  
\usepackage{amsmath}
\usepackage{amssymb}

\newtheorem{theorem}{Theorem}
\newtheorem{corollary}{Corollary}
\newtheorem{proposition}{Proposition}

\newenvironment{definition}
{\smallskip\noindent{\bf Definition\/}:}{\smallskip\par}

\newenvironment{example}
{\smallskip\noindent{\bf Example\/}.}{\smallskip\par}
\newenvironment{remark}
{\smallskip\noindent{\bf Remark\/}.}{\smallskip\par}
\newenvironment{remarks}
{\smallskip\noindent{\bf Remarks\/}.}{\smallskip\par}
\newenvironment{proof}
{\noindent{\it Proof\/}.}{{ \hfill $\Box$}\smallskip\par}
\newenvironment{Proof}
{\noindent{\it Proof\/}}{{ \hfill $\Box$}\smallskip\par}

\newcommand{\CC}{{\Bbb C}}

\newcommand{\RR}{{\Bbb R}}

\newcommand{\eps}{\varepsilon}

\title{Radial index and Euler obstruction of a 1-form on a singular variety}
\author{W.~Ebeling and S.~M.~Gusein-Zade
\thanks{Partially supported by the DFG-programme ''Global methods in
complex geometry'' (Eb 102/4--2), grants RFBR--01--01--00739, NSh--1972.2003.1.
}
}

\date{}

\begin{document}

\maketitle

\begin{abstract} A notion of the radial index of an isolated singular point of a 1-form on a
singular (real or complex) variety is discussed. For the differential of a function it is
related to the Euler characteristic of the Milnor fibre of the function. A connection between
the radial index and the local Euler obstruction of a 1-form is described. 
This gives an expression for the local Euler
obstruction of the differential of a function in terms of Euler characteristics of some
Milnor fibres.
\end{abstract}

\section*{Introduction}

An isolated singular point (zero) of a vector field or of a 1-form on a smooth manifold (real
or complex) has a natural integer invariant -- the index. There are several generalizations
of this notion: an index of a vector field on an isolated singularity \cite{ASV, EG1}, an index
of a vector field \cite{GSV, SS} or of a 1-form \cite{EG2} on an isolated complete intersection
singularity, \dots.

Let $(X, 0)\subset(\RR^N, 0)$ be the germ of a purely $n$-dimensional analytic variety
and let $X = \bigcup V_i$ be a Whitney stratification of it. Let
$v$ be a continuous vector field on a neighbourhood of the origin in $\RR^N$ tangent to $X$. The
last fact means that the vector field $v$ is tangent to each stratum $V_i$ of the Whitney
stratification of
$(X, 0)$. Suppose that $v(x)\ne 0$ for all points $x$ from a punctured neighbourhood of the
origin in $X$. In \cite{EG1} there was defined a notion of an index of the vector field $v$ on
$X$ at the origin which we shall call the {\em radial index}. (For a germ $(X, 0)$ with an
isolated singular point at the origin this notion was defined and discussed in
\cite{ASV, EG1}.) 

The condition for a vector field to be tangent to a germ of a
singular variety is a very restrictive one. For example, holomorphic vector fields with
isolated zeroes exist on a complex analytic variety with an isolated singularity
\cite{BG}, but not in general. 

\begin{example}  
Consider the surface
$X$ in $\CC^3$ given by the equation 
$$xy(x-y)(x+zy)=0.$$
It has singularities on the line $\{x=y=0\}$. It
is a family of four lines in the $(x,y)$-plane with different cross ratios. Then any
holomorphic vector field tangent to $X$ vanishes on the line $\{x=y=0\}$, because the
translation along such a vector field has to preserve the cross ratio of the lines.
\end{example}

However, 1-forms with isolated singular points always exist. Therefore
(and because of some other advantages: see the discussion in
\cite{EG2}) we prefer to consider 1-forms instead of vector fields. 

Suppose that $\omega$ is a continuous 1-form defined in a neighbourhood of the
origin in
$\RR^N$. A singular point of the 1-form $\omega$ on $X$ is a singular point of its
restriction to a stratum of the Whitney stratification of $X$. (One should consider all points
of a zero-dimensional manifold as singular ones.) Suppose that the 1-form $\omega$ 
has an isolated singular point on $X$ at the origin. In this situation we define the
notion of the radial index of the 1-form $\omega$ on $(X, 0)$ (or rather adapt one for vector
fields from
\cite{EG1}). We also define a similar concept for complex 1-forms on a germ of a
complex analytic subvariety of $(\CC^N, 0)$.

For a germ $f$ of an analytic function on the germ $(X, 0)$ of an analytic variety, there
exists the notion of the Milnor fibration and of the Milnor fibre (or of the Milnor fibres in
the real setting: the positive and the negative ones): see, e.g., \cite{Le, Var}. If the
function $f$ has an isolated critical point on $X$ at the origin (i.e. if its differential
$df$ has an isolated singular point there) we describe a connection between the radial index
of the differential
$df$ and the Euler characteristic of the Milnor fibre of the germ $f$. In the described
situation there is also an appropriately adapted notion of the local Euler obstruction of a
1-form on the germ of a complex analytic variety \cite{BMPS}. We describe a connection between
the radial index and the local Euler obstruction. This gives an expression for the local Euler
obstruction of (the differential of) a function in terms of Euler characteristics of some
Milnor fibres. 

After the paper was put to the Math.~ArXiv, J.~Sch\"urmann informed us that our definition of
the radial index can be interpreted as a microlocal intersection number defined in \cite{KS}.
With this definition part of the statements of this paper (essentially for $\omega=df$ where
$f$ is an analytic function) can be deduced from general statements of \cite{KS} (see also
\cite{Sch1, Sch2}).

\section{The radial index of a 1-form}

First let $(X, 0)\subset(\RR^N, 0)$ be the germ of a real analytic variety and let $\omega$ be a
(continuous) 1-form on a neighbourhood of the origin in $\RR^N$. 

\begin{definition} The 1-form $\omega$ is {\em radial} on $(X, 0)$ if, for an arbitrary
nontrivial analytic arc $\varphi: (\RR, 0)\to (X, 0)$ on $(X, 0)$, the value of the 1-form
$\omega$ on the tangent vector $\dot{\varphi}(t)$ is positive for positive $t$ (small enough).
\end{definition}

Radial 1-forms exist: e.g.\ the differential $d\phi$ for the germ of an analytic function
$\phi: (\RR^n,0) \to (\RR,0)$ with a strict local minimum at the origin (one can take
$\phi=r^2$  where $r$ is the distance from the origin in $\RR^N$).

Let $X = \bigcup_{i=0}^q V_i$ be a Whitney stratification of the germ $(X, 0)$,
$V_0=\{0\}$. Let $\omega$ be an (arbitrary continuous) 1-form on a neighbourhood of the origin
in
$\RR^N$ with an isolated singular point on $(X, 0)$ at the origin. Let $\eps>0$ be small enough
so that in the closed ball $B_\eps$ of radius $\eps$ centred at the origin in $\RR^N$ the
1-form $\omega$ has no singular points on $X\setminus\{0\}$. It is easy to see that there
exists a 1-form $\widetilde\omega$ on $\RR^N$ such that:
\begin{itemize}
\item[1)] The 1-form
$\widetilde\omega$ coincides with the 1-form $\omega$ on a neighbourhood of the sphere $S_\eps
= \partial B_\eps$. 
\item[2)] The 1-form $\widetilde\omega$ is radial on $(X, 0)$ at the origin. 
\item[3)] In a
neighbourhood of each singular point $x_0\in (X\cap B_\eps)\setminus\{0\}$, $x_0 \in
V_i$, $\dim V_i=k$, the 1-form $\widetilde{\omega}$ looks as follows. There exists a
(local) analytic diffeomorphism $h: (\RR^N, \RR^k,0) \to (\RR^N,V_i, x_0)$ such that $h^\ast
\widetilde{\omega} = \pi_1^\ast \widetilde{\omega}_1 + \pi_2^\ast \widetilde{\omega}_2$, where
$\pi_1$ and $\pi_2$ are the natural projections $\pi_1: \RR^N \to \RR^k$
and $\pi_2: \RR^N \to \RR^{N-k}$ respectively, 
$\widetilde{\omega}_1$ is the germ of a 1-form on $(\RR^k,0)$ with an isolated singular point
at the origin, and $\widetilde{\omega}_2$ is a radial 1-form on $(\RR^{N-k},0)$.
\end{itemize} 

\begin{remark} One can demand that the 1-form $\widetilde\omega_1$ has a non-degenerate
singular point (and therefore ${\rm ind}_0 \, \widetilde\omega_1 = \pm 1$), however, this is
not necessary for the definition.
\end{remark}


\begin{definition} The {\em radial index} ${\rm ind}_{X,0} \, \omega$ of the 1-form $\omega$ on
the variety
$X$ at the origin is the sum
$$ 1 + \sum_{i=1}^q\sum_{Q\in {\rm Sing}\, \widetilde{\omega} \cap V_i} {\rm ind}_Q \,
\widetilde\omega_{\vert V_i}
$$
where the sum is taken over all singular points of the 1-form $\widetilde{\omega}$ on $(X
\setminus \{0\}) \cap B_\eps$.
\end{definition}

\begin{remarks} 1. One can write the definition as 
$$
\sum_{i=0}^q\sum_{Q\in {\rm Sing}\, \widetilde{\omega} \cap V_i} {\rm ind}_Q \,
\widetilde\omega_{\vert V_i}
$$ 
(i.e.\ over all singular points of the 1-form $\widetilde{\omega}$ on $X \cap B_\eps$) if one
assumes that the index of a 1-form on a point (a zero-dimensional manifold) is equal to 1. 

2. Since no other indices will be considered in this paper, we shall usually omit the word
{\em radial}.
\end{remarks}

To prove that this notion is well defined, we shall show that not only the total sum 
$$
\sum_{i=0}^q\sum_{Q\in {\rm Sing}\, \widetilde{\omega} \cap V_i} {\rm ind}_Q \,
\widetilde\omega_{\vert V_i}
$$ 
for all the strata is well defined but also the sum
$$
\sum_{Q\in {\rm Sing}\, \widetilde{\omega} \cap V_i} {\rm ind}_Q \, \widetilde\omega_{\vert
V_i}$$ for each single stratum. 

To prove this for a stratum $V_i$ we suppose that two 1-forms $\widetilde{\omega}$ and
$\widetilde{\omega}'$ have singular points of the described type. One can change each of
these forms in a neighbourhood of $\overline{V_i} \setminus V_i$ by such a continuous
deformation that in the course of it no new singular points arise  and the value
of the new forms on a fixed  outward looking normal vector field on the boundary of a tubular
neighbourhood of 
$\overline{V_i} \setminus V_i$ is positive. This can be achieved locally by adding 
a 1-form radial in the normal direction to a stratum (say, the
differential of the squared distance from the stratum defined in local coordinates). After
that, the new 1-forms $\widetilde{\omega}$ and $\widetilde{\omega}'$, as non-vanishing
sections of the cotangent bundle, are homotopic to each other in a neighbourhood of the
boundary of the stratum $V_i \cap B_\eps$. This implies that the sums of the indices of the
singular points of the 1-forms
$\widetilde{\omega}$ and $\widetilde{\omega}'$ inside $V_i \cap B_\eps$ coincide.

Just because of the definition, the (radial) index satisfies the law of conservation of
number, i.e., if a 1-form $\omega'$ with isolated singular points on $X$ is close to the 1-form
$\omega$, then
$${\rm ind}_{X,0} \, \omega = \sum_{Q \in {\rm Sing}\, \omega'} {\rm ind}_{X,Q}\,  \omega'$$
where the sum on the right hand side is over all singular points $Q$ of the 1-form $\omega'$ on
$X$ in a neighbourhood of the origin. This implies the following statement. Let $X$ be a
compact analytic subset of $\RR^N$ with a Whitney stratification. Let $\omega$ and $\omega'$
be 1-forms on $\RR^N$ with isolated singular points on $X$.

\begin{proposition} \label{prop1}
One has
$$\sum_Q {\rm ind}_{X,Q}\,  \omega = \sum_{Q'} {\rm ind}_{X,Q'}\,  \omega'$$
where the sums are taken over all singular points of the 1-forms $\omega$ and $\omega'$
respectively. 
\end{proposition}

\begin{proof}
This follows from the law of conservation of number applied to a family of
1-forms (with isolated singularities) connecting the 1-forms $\omega$ and $\omega'$.
\end{proof}

In fact one has a more precise statement:

\begin{theorem}[Poincar\'e-Hopf] One has
$$\sum_Q {\rm ind}_{X,Q}\,  \omega = \chi(X)$$
where $\chi(X)$ denotes the Euler characteristic of the set $X$.
\end{theorem}

We postpone the proof to Section~\ref{sect2}.

\smallskip

Now let $(X,0) \subset (\CC^N,0)$ be the germ of a complex analytic variety of pure dimension
$n$ and let $\omega$ be a (complex and, generally speaking, continuous) 1-form on a
neighbourhood of the origin in $\CC^N$. In fact there is a one-to-one correspondence between
complex 1-forms on a complex manifold $M^n$ (say, on $\CC^N$) and real 1-forms on it
(considered as a real $2n$-dimensional manifold). Namely, to a complex 1-form $\omega$ one
associates the real 1-form $\eta= {\rm Re}\, \omega$; the 1-form $\omega$ can be restored from
$\eta$ by the formula $\omega(v)=\eta(v)-i\eta(iv)$ for $v \in T_x M^n$. This means that the
index of the real 1-form ${\rm Re}\, \omega$ is an invariant of the complex 1-form
$\omega$ itself. However, on a smooth manifold ${\rm ind}_{M^n,x}\,  {\rm Re}\, \omega$
does not coincide with the usual index of the singular point $x$ of the 1-form $\omega$, but
differs from it by the coefficient $(-1)^n$. (E.g., the index of the (complex analytic)
1-form $\omega =
\sum_{j=1}^n x_j dx_j$ ($(x_1,\ldots, x_n)$ being the coordinates of $\CC^n$) is equal to 1,
whence the index of the real 1-form ${\rm Re}\, \omega = \sum_{j=1}^n u_j du_j -
\sum_{j=1}^n v_j dv_j$ ($x_j=u_j+iv_j$) is equal to $(-1)^n$.) This explains the following
definition.

\begin{definition}
The {\em (complex radial) index} ${\rm ind}^{\CC}_{X,0}\,  \omega$ of the complex 1-form
$\omega$ on
$X$ at the origin is $(-1)^n$ times the index of the real 1-form ${\rm Re}\, \omega$ on $X$:
$${\rm ind}^{\CC}_{X,0}\,  \omega = (-1)^n \,{\rm ind}_{X,0}\,  {\rm Re}\, \omega.$$
\end{definition}

When it is clear from the context, the upper index $\CC$ will be omitted.

\begin{example} Let the variety $(X,0)$ have an isolated singularity at the origin and let
$\eps>0$ be small enough so that the 1-form $\omega$ has no singular points on the
intersection of $X$ with the ball $B_\eps$ of radius $\eps$ centred at the origin in $\CC^N$
outside of the origin. Let $\pi : (\widetilde{X}, D)
\to (X,0)$ be a resolution of the singularity $(X,0)$ which is an isomorphism outside of the
origin. The lifting $\pi^\ast\omega$ of the 1-form $\omega$ to the space $\widetilde{X}$ of
the resolution has no zeroes outside of the exceptional divisor $D$. Let ${\rm Obst}_\pi \,
\omega$ be the obstruction to extend the non-zero section $\pi^\ast \omega$ of the
cotangent bundle of $\widetilde{X}$ from a neighbourhood of the preimage $\pi^{-1}(S_\eps
\cap X)$  ($S_\eps= \partial B_\eps$) to the whole preimage of $B_\eps \cap X$. Since the
obstruction  ${\rm Obst}_\pi \, \omega$ (for a fixed resolution $\pi$) satisfies the law of
conservation of number as well and for the trivial resolution of a smooth manifold  
${\rm Obst}_\pi \, \omega$ coincides with the index ${\rm ind}_{X,0} \, \omega$, the
difference 
$${\rm Obst}_\pi \, \omega - {\rm ind}_{X,0} \, \omega$$
does not depend on the 1-form $\omega$. Moreover, the manifold $\widetilde{X}$ can be
considered as a smoothing of the variety $X$ in the sense of \cite{EG1,EG2}. Therefore
$${\rm Obst}_\pi \, \omega - {\rm ind}_{X,0} \, \omega = (-1)^n (\chi(D)-1)$$
(cf.\ \cite[Proposition~2]{EG2}; note that the sign $(-1)^n$ appears because we consider here
the complex analytic situation).
\end{example}

In the sequel it will convenient to denote $\chi(Z)-1$ by $\overline{\chi}(Z)$ and to call it
the {\em reduced} (modulo a point) Euler characteristic of the topological space $Z$ (though,
strictly speaking, this name is only correct for a non-empty space $Z$).

\section{The Milnor fibre of a function and the index of its differential} \label{sect2}
We shall use the following two statements of stratified Morse theory. The first one is
very general. We formulate the second one in the simplest form necessary for us. The
formulations here are made in the global setting. The necessary changes for the local case
(i.e.\ inside a ball) are obvious.

Let $X$ be a compact analytic subset of $\RR^N$ with a Whitney stratification and let $f:
\RR^N \to \RR$ be a smooth function with isolated critical points on $X$. 
For a subset $A \subset \RR$, 
let $M_A =
\{ x \in X \, | \, f(x) \in A \}$. 

\begin{proposition} \label{prop2} 
Suppose that the function $f$ has no critical values in the segment
$[a,b]$. Then $M_{\{a\}}$ is a homotopy retract of $M_{[a,b]}$.
\end{proposition}

See the proof in \cite{GM}. It is a corollary of Thom's first isotopy lemma.

\smallskip

Suppose that there is only one critical point $x_0$ of the function $f$ on $X$ with the
critical value $c$ and that in a neighbourhood of the point $x_0$ the function $f$ has the
following form. Let $V_i$ be the stratum of $X$ containing $x_0$ and let $k$ be the dimension
of $V_i$. Then there exists a local analytic diffeomorphism $h: (\RR^N,\RR^k,0) \to (\RR^N,
V_i, x_0)$. Let $y_1, \ldots , y_N$ be the coordinates in $\RR^N$, $\RR^k = \{ y \in \RR^N \,
| \, y_{k+1} = \ldots = y_N=0\}$. Suppose that the function $h^\ast f$ ($=f \circ h$) is
equal to the sum
$f_1(y_1, \ldots , y_k) + f_2(y_{k+1}, \ldots , y_N)$
where 
the function $f_1 : (\RR^k,0) \to (\RR,0)$ has a non-degenerate critical point at the origin
with Morse index $m$ and  $f_2 : (\RR^{N-k},0) \to (\RR,0)$ is an analytic germ with a strict
local minimum at the origin (e.g. $f_2(y_{k+1}, \ldots , y_N)= y_{k+1}^2+ \ldots + y_N^2$).

\begin{proposition} \label{prop3}
For $\eps >0$ small enough, the space $M_{[c - \eps, c + \eps]}$ can be
homotopically retracted to the subspace $M_{\{c-\eps\}}$ with one cell of dimension $m$
attached.
\end{proposition}

The proof is very easy. See \cite{GM} for similar, but stronger statements.

\smallskip

\begin{Proof}{\it of the Poincar\'e-Hopf theorem.} Since the sum $\sum {\rm ind}_{X,Q} \,
\omega$ does not depend on the 1-form $\omega$, one can take $\omega = df$, where the
smooth function
$f$ has only critical points of the type described before Proposition~\ref{prop3}. Then
Propositions~\ref{prop2} and \ref{prop3} yield the statement.
\end{Proof}

Let $(X,0) \subset (\RR^N,0)$ be the germ of a real analytic variety of pure dimension $n$ and
let $f:(X,0) \to (\RR,0)$ be the germ of an analytic function on $X$. The germ $f$ is the
restriction to $X$ of the germ of an analytic function on $(\RR^N,0)$ (also denoted by $f$).
Suppose that $f$ has an isolated critical point on $X$ at the origin, i.e.\ its differential
$df$ has an isolated singular point there.
For a small enough positive $\eps$, the restriction of the function $f$ to the intersection $X
\cap B_\eps$ defines a map $f: X \cap B_\eps \to \RR$ which is a locally trivial fibration
over a punctured neighbourhood of zero in $\RR$, generally speaking with different fibres over
the positive and negative parts (see, e.g., \cite{Var}). Let $M_f^+$ and $M_f^-$ be the fibres
of this fibration which can be called positive and negative Milnor fibres of $f$ respectively.

\begin{theorem} \label{theo2} {\rm (cf.\ \cite[Example~2.6]{EGS})}
One has
$${\rm ind}_{X,0} \, df = -\overline{\chi}(M_f^-).$$
\end{theorem}

\begin{proof} Let $\delta$ be positive and small enough ($0< \delta \ll \eps$) so that
$M_{\{-\delta\}} = f^{-1}(\delta) \cap B_\eps$ is the negative Milnor fibre of the function $f$
and
$M_{[-\delta,
\delta]}=f^{-1}([-\delta, \delta]) \cap B_\eps$ is contractible. One can perturb the function
$f$ in a small neighbourhood of the origin (in which the absolute value of the function $f$ is
smaller than $\delta$) so that it will only have singular points of the type described before
Proposition~\ref{prop3}. Now the local versions of Propositions~\ref{prop2} and \ref{prop3}
yield the statement.
\end{proof}

Let $(X,0) \subset (\CC^N,0)$ be the germ of a complex analytic variety of pure dimension $n$
and let $f: (X,0) \to (\CC,0)$ be the germ of a holomorphic function with an isolated critical
point at the origin.

\begin{theorem} \label{theo3}
One has
$${\rm ind}^{\CC}_{X,0} \, df = (-1)^{n-1}\overline{\chi}(M_f).$$
\end{theorem}

\begin{proof} By \cite[Proposition~2.A.3]{GM}, the Milnor fibre $M^\pm_{{\rm Re}\, f}$ of the
the real part ${\rm Re}\, f$ of the function $f$ on
$(X,0)$ is homeomorphic to the cartesian product $M_f \times I$ of the Milnor
fibre $M_f$ of the function $f$ and the segment $I=[0,1]$.
\end{proof}

\section{The local Euler obstruction of a 1-form}
In \cite{BMPS}, there is introduced the notion of the local Euler obstruction ${\rm
Eu}_{f,X}(0)$ of a holomorphic function $f: (X,0) \to (\CC,0)$ with an isolated critical point
on the germ of a complex analytic variety $(X,0)$. It is defined through an appropriately
constructed gradient vector field tangent to the variety $X$. We adapt the definition to the
case of a 1-form.

Let $(X,0) \subset (\CC^N,0)$ be the germ of a complex analytic variety with a Whitney
stratification
$X=\bigcup_{i=0}^q V_i$, $V_0=\{0\}$, and let $\omega$ be a 1-form on a neighbourhood of the
origin in $\CC^N$ with an isolated singular point on $X$ at the origin. Let $\eps>0$ be small
enough such that the 1-form
$\omega$ has no singular points on $X \setminus \{0\}$ inside the ball $B_\eps$. 
Let $\nu : \widehat{X} \to X$ be the Nash
transformation of the variety $X$ defined as follows. Let $G(n,N)$ be the Grassmann manifold
of $n$-dimensional vector subspaces of $\CC^N$. For a suitable neighbourhood $U$ of the origin
in $\CC^N$, there is a natural map $\sigma : X_{\rm reg} \cap U \to U \times G(n,N)$ which
sends a point $x$ to $(x, T_x X_{\rm reg})$ ($X_{\rm reg}$ is the non-singular part of $X$).
The Nash transform $\widehat{X}$ is the closure of the image ${\rm Im}\, \sigma$ of the map
$\sigma$ in $U
\times G(n,N)$. The Nash bundle $\widehat{T}$ over $\widehat{X}$ is a vector bundle of rank $n$
which is the pullback of the tautological bundle on the Grassmann manifold $G(n,N)$. There is a
natural lifting of the Nash transformation to a bundle map from the Nash bundle
$\widehat{T}$ to the restriction of the tangent bundle $T\CC^N$ of $\CC^N$ to $X$. This is
an isomorphism of
$\widehat{T}$ and $TX_{\rm reg} \subset T\CC^N$ over the regular part $X_{\rm reg}$ of $X$.
The 1-form $\omega$ gives rise to a section $\widehat{\omega}$ of the dual Nash bundle
$\widehat{T}^\ast$ over the Nash transform $\widehat{X}$ without zeroes outside of the
preimage of the origin.

\begin{definition} The {\em local Euler obstruction} ${\rm Eu}_{X,0} \, \omega$ of the 1-form
$\omega$ on $X$ at the origin is the obstruction to extend the non-zero section
$\widehat{\omega}$ from the preimage of a neighbourhood of the sphere $S_\eps= \partial
B_\eps$ to the preimage of its interior, more precisely its value (as an element of
$H^{2n}(\nu^{-1}(X\cap B_\eps), \nu^{-1}(X \cap S_\eps))$\,) on the fundamental class of the
pair $(\nu^{-1}(X\cap B_\eps), \nu^{-1}(X \cap S_\eps))$.
\end{definition}

The word {\em local} will usually be omitted.

\begin{remark} The local Euler obstruction can also be defined for a real 1-form on the germ
of a real analytic variety if the last one is orientable in an appropriate sense.
\end{remark}

\begin{example} Let $\omega = df$ for the germ $f$ of a holomorphic function on $(\CC^N,0)$.
Then 
${\rm Eu}_{X,0} \, df$ differs from the Euler obstruction ${\rm Eu}_{f,X}(0)$ defined in
\cite{BMPS} by the sign $(-1)^n$. The reason is that for the germ of a holomorphic function
with an isolated critical point on $(\CC^n,0)$ one has ${\rm Eu}_{f,X}(0) = (-1)^n \mu_f$ (see
\cite[Remark~3.4]{BMPS}), whence ${\rm Eu}_{X,0} \, df= \mu_f$ ($\mu_f$ is the Milnor number
of the germ $f$). E.g., for $f(x_1, \ldots , x_n)=x_1^2 + \ldots + x_n^2$ the obstruction  
${\rm Eu}_{f,X}(0)$ is the index of the vector field $\sum_{i=1}^n \overline{x}_i \partial/
\partial x_i$ (which is equal to $(-1)^n$), but the obstruction ${\rm Eu}_{X,0} \, df$ is the
index of the (holomorphic) 1-form $\sum_{i=1}^n x_i dx_i$ which is equal to 1.
\end{example}

The Euler obstruction of a 1-form satisfies the law of conservation of number (just as
the radial index). Moreover, on a smooth variety the Euler obstruction and the radial
index coincide. This implies the following statement (cf.\ \cite[Theorem~3.1]{BMPS}).

\begin{proposition} \label{prop4}
Let $(X,0) \subset (\CC^N,0)$ have an isolated singular point at the origin and let $\ell:
\CC^N
\to
\CC$ be a generic linear function. Then
$${\rm ind}_{X,0} \, \omega - {\rm Eu}_{X,0} \, \omega = {\rm ind}_{X,0} \, d\ell =
(-1)^{n-1}\overline{\chi}(M_\ell),$$
where $M_\ell$ is the Milnor fibre of the linear function $\ell$ on $X$.
In particular
$${\rm Eu}_{X,0} \, df = (-1)^n (\chi(M_\ell) - \chi(M_f)).$$
\end{proposition}

\begin{proof} One can deform the 1-form $\omega$ to a 1-form $\widetilde{\omega}$ such that
$\ell=\widetilde{\omega}(0)$ is a generic linear function on $\CC^N$. 
Since at the singular points of $\widetilde{\omega}$ on $X$ outside of the
origin the Euler obstructions and the radial indices coincide, one has
$${\rm ind}_{X,0} \, \omega - {\rm Eu}_{X,0} \, \omega = {\rm ind}_{X,0} \, \widetilde{\omega}
- {\rm Eu}_{X,0} \, \widetilde{\omega}.$$
It is easy to see that both the Euler obstruction and the index of the 1-form
$\widetilde{\omega}$ on $X$ at the origin coincide with those of the differential of the linear
function
$\ell$. The Euler obstruction ${\rm Eu}_{X,0} \, d\ell$ is equal to zero (see
\cite[Lemma~1.3]{BLS}) and 
$${\rm ind}_{X,0} \, d\ell = (-1)^{n-1} \overline{\chi}(M_\ell)$$
by Theorem~\ref{theo3}. This yields the statement.
\end{proof}

Now let $(X,0) \subset (\CC^N,0)$ be an arbitrary germ of an analytic variety with a Whitney
stratification $X=\bigcup_{i=0}^q V_i$, $V_0=\{0\}$. For a stratum $V_i$, $i=0, \ldots , q$,
let $N_i$ be the normal slice in the variety $X$ to the stratum $V_i$ ($\dim N_i =
\dim X- \dim V_i$) at a point of the stratum $V_i$ and let $n_i$ be the index of a generic
(non-vanishing) 1-form $d\ell$ on $N_i$:
$$n_i = (-1)^{\dim N_i-1} \overline{\chi}(M_{\ell|_{N_i}}).$$
In particular for an open stratum $V_i$ of $X$, $N_i$ is a point and $n_i=1$.
The strata $V_i$ of $X$ are partially ordered: $V_i \prec
V_j$ (we shall write $i \prec j$) iff $V_i \subset \overline{V_j}$ and $V_i \neq V_j$; $i
\preceq j$ iff $i \prec j$ or $i=j$. For two strata $V_i$ and $V_j$ with $i \preceq j$, let
$N_{ij}$ be the normal slice of the variety $\overline{V_j}$ to the stratum $V_i$ at a point
of it ($\dim N_{ij}= \dim V_j - \dim V_i$, $N_{ii}$ is a point) and let $n_{ij}$ be the index
of a generic (non-vanishing) 1-form $d\ell$ on $N_{ij}$: $n_{ij}=(-1)^{\dim
N_{ij}-1}\, \overline{\chi}(M_{\ell\vert_{N_{ij}}})$,
$n_{ii}=1$.  Let us define the Euler obstruction ${\rm
Eu}_{Y,0} \, \omega$ to be equal to 1 for a zero-dimensional variety
$Y$ (in particular ${\rm Eu}_{\overline{V_0},0}\, \omega =1$, ${\rm Eu}_{N_{ii},0}\,
\omega =1$).

\begin{theorem} \label{theo4}
One has
$${\rm ind}_{X,0} \, \omega = \sum_{i=0}^q n_i \cdot {\rm Eu}_{\overline{V_i},0}\, \omega.$$
\end{theorem}

\begin{proof} There exists a 1-form $\widetilde{\omega}$ which coincides with $\omega$ in a
neighbourhood of the sphere $S_\eps$ such that in a neighbourhood of each singular
point $x_0$, $x_0 \in V_i$, $\dim V_i=k$, it looks as follows. There exists a local
biholomorphism $h: (\CC^N, \CC^k,0) \to (\CC^N, V_i, x_0)$ and, in coordinates $y_1, \ldots ,
y_N$ in $\CC^N$ with $\CC^k=\{y_{k+1}= \ldots = y_N=0\}$, one has $h^\ast \widetilde{\omega} =
\pi_1^\ast \omega_1 + d\ell$ where $\omega_1$ is a 1-form on $(\CC^k,0)$ with a non-degenerate
singular point at the origin (and therefore ${\rm ind}_{\CC^k,0}\, \omega_1= \pm 1$), $\pi_1$
is the projection $\pi_1: \CC^N \to \CC^k$, and $\ell=\ell(y_{k+1}, \ldots , y_N)$ is a generic
linear function on $\CC^{N-k}$. For $i \preceq j$ one  obviously has
\begin{eqnarray*}
{\rm ind}_{\overline{V_j},x_0} \, \widetilde{\omega} & = & {\rm ind}_{\CC^k,0}\, \omega_1 \cdot
{\rm ind}_{N_{ij},x_0}\, d\ell,\\
{\rm Eu}_{\overline{V_j},x_0} \, \widetilde{\omega} & = & {\rm ind}_{\CC^k,0}\, \omega_1 \cdot
{\rm Eu}_{N_{ij},x_0}\, d\ell.
\end{eqnarray*}
Since ${\rm Eu}_{N_{ij},x_0}\, d\ell=0$ for $i \prec j$ \cite[Lemma~1.3]{BLS} and ${\rm
Eu}_{N_{ii},x_0}\, d\ell = {\rm ind}_{N_{ii},x_0}\, d\ell=1$, one has
$$\sum_{x_0 \in V_j} {\rm ind}_{V_j,x_0} \, \widetilde{\omega}|_{V_j} = 
{\rm Eu}_{\overline{V_j},0} \, \omega,$$
$$
{\rm ind}_{X,0}\, \omega  =  \sum_{j=0}^q \sum_{x_0 \in V_j} {\rm ind}_{X,x_0}\,
\widetilde{\omega} = \sum_{j=0}^q \sum_{x_0 \in V_j} n_j \cdot {\rm ind}_{V_j,x_0}\,
\widetilde{\omega}|_{V_j} = \sum_{j=0}^q n_j \cdot {\rm Eu}_{\overline{V_j},0} \, \omega.
$$
\end{proof}

To write an "inverse" of the formula of Theorem~\ref{theo4}, suppose that the variety $X$ is
irreducible and $X=\overline{V_q}$. (Otherwise one can permit $X$ to be reducible, but also
permit the open stratum $V_q$ to be not connected and dense; this does not change anything in
Theorem~\ref{theo4}.)  Let
$m_{ij}$ be the (M\"obius) inverse of the function
$n_{ij}$ on the partially ordered set of strata, i.e.
$$\sum_{i \preceq j \preceq k} n_{ij}m_{jk} = \delta_{ik}$$
(see, e.g., \cite{Hall}). For $i \prec j$ one has
\begin{eqnarray*}
m_{ij} &  = &  \sum_{i=k_0 \prec k_1 \prec \ldots \prec k_r =j} (-1)^{r}
n_{k_0k_1}n_{k_1k_2} \ldots n_{k_{r-1}k_r}\\
& = & 
(-1)^{\dim X -\dim V_i} \sum_{i=k_0 \prec \ldots \prec k_r =j}
\overline{\chi}(M_{\ell\vert_{N_{k_0k_1}}}) \cdot \ldots \cdot
\overline{\chi}(M_{\ell\vert_{N_{k_{r-1}k_r}}})\ .
\end{eqnarray*}

\begin{corollary} 
One has
$${\rm Eu}_{X,0}\, \omega = \sum_{i=0}^q m_{iq} \cdot {\rm ind}_{\overline{V_i},0} \, \omega.$$
In particular 
\begin{eqnarray*}
\lefteqn{{\rm Eu}_{X,0}\, df = (-1)^{\dim X -1} \times}\\
&& \left( \overline{\chi}(M_{f\vert_X}) + \sum_{i=0}^{q-1}
\overline{\chi}(M_{f\vert_{\overline{V_i}}})  
\sum_{i=k_0 \prec \ldots \prec k_r =q}
\overline{\chi}(M_{\ell\vert_{N_{k_0k_1}}})  \ldots 
\overline{\chi}(M_{\ell\vert_{N_{k_{r-1}k_r}}}) \right).
\end{eqnarray*}
\end{corollary}

\bigskip
\noindent Universit\"{a}t Hannover, Institut f\"{u}r Mathematik \\
Postfach 6009, D-30060 Hannover, Germany \\
E-mail: ebeling@math.uni-hannover.de\\

\medskip
\noindent Moscow State University, Faculty of Mechanics and Mathematics\\
Moscow, 119992, Russia\\
E-mail: sabir@mccme.ru

\end{document}